\documentclass{conm-p-l}

\newtheorem{theorem}{Theorem}[section]
\newtheorem{lemma}[theorem]{Lemma}

\theoremstyle{definition}
\newtheorem{definition}[theorem]{Definition}

\theoremstyle{remark}
\newtheorem{remark}[theorem]{Remark}

\numberwithin{equation}{section}

\newcommand{\LL}{{\mathcal L}}
\newcommand{\tla}{\tilde{\lambda}}

\newcommand{\hk}{\hat{k}}
\newcommand{\Z}{Z^2/\{0\}}
\newcommand{\E}{{\mathcal E}}
\newcommand{\HH}{{\mathcal H}}

\newcommand{\e}{\varepsilon}
\newcommand{\vth}{\vartheta}
\newcommand{\U}{{\mathcal U}}

\renewcommand{\k}{\kappa}
\newcommand{\ga}{\gamma}
\newcommand{\Ga}{\Gamma}

\newcommand{\dl}{\delta}
\newcommand{\Dl}{\Delta}
\renewcommand{\th}{\theta}
\newcommand{\ra}{\rightarrow}
\newcommand{\al}{\alpha}
\newcommand{\be}{\beta}
\newcommand{\sg}{\sigma}
\newcommand{\Sg}{\Sigma}

\newcommand{\pa}{\partial}
\newcommand{\z}{\zeta}

\newcommand{\La}{\Lambda}

\newcommand{\la}{\lambda}
\newcommand{\bq}{\bar{q}}

\newcommand{\nid}{\noindent}

\newcommand{\om}{\omega}
\newcommand{\Om}{\Omega}
\newcommand{\na}{\nabla}

\newcommand{\W}{{\mathcal W}}

\def\maprightu#1{\smash{
    \mathop{\longrightarrow}\limits^{#1}}}
\def\maprightd#1{\smash{
    \mathop{\longrightarrow}\limits_{#1}}}
\def\mapdownl#1{
    \llap{$\vcenter{\hbox{$\scriptstyle#1$}}$}\Big\downarrow}
\def\mapdownr#1{\Big\downarrow
    \rlap{$\vcenter{\hbox{$\scriptstyle#1$}}$}}

\begin{document}

\title{Chaos in Partial Differential Equations}

\author{Yanguang (Charles)  Li}
\address{Department of Mathematics, University of Missouri, 
Columbia, MO 65211}
\curraddr{}
\email{cli@math.missouri.edu}
\thanks{This work was partially supported by a Guggenheim Fellowship}

%\author{Author Two}
%\address{Mathematical Research Section, School of Mathematical Sciences,
%Australian National University, Canberra ACT 2601, Australia}

%    General info
\date{}



\begin{abstract}
This is a survey on Chaos in Partial Differential Equations. First 
we classify soliton equations into three categories: 1. (1+1)-dimensional
soliton equations, 2. soliton lattices, 3. (1+n)-dimensional
soliton equations ($n\geq 2$). A systematic program has been established 
by the author and collaborators, for proving the existence
of chaos in soliton equations under perturbations. For each category, 
we pick a representative to present the results. Then we review some 
initial results on 2D Euler equation.
\end{abstract}

\maketitle

%\section*{}
%This is an example of an unnumbered first-level heading.

%\specialsection*{This is a Special Section Head}
%This is an example of a special section head%
%%%%%%%%%%%%%%%%%%%%%%%%%%%%%%%%%%%%%%%%%%%%%%%%%%%%%%%%%%%%%%%%%%%%%%%%
%\footnote{Here is an example of a footnote. Notice that this footnote
%text is running on so that it can stand as an example of how a footnote
%with separate paragraphs should be written.
%\par
%And here is the beginning of the second paragraph.}%
%%%%%%%%%%%%%%%%%%%%%%%%%%%%%%%%%%%%%%%%%%%%%%%%%%%%%%%%%%%%%%%%%%%%%%%%

%\section{This is a numbered first-level section head}
%This is an example of a numbered first-level heading.

%\subsection{This is a numbered second-level section head}
%This is an example of a numbered second-level heading.

%\subsection*{This is an unnumbered second-level section head}
%This is an example of an unnumbered second-level heading.

%\subsubsection{This is a numbered third-level section head}
%This is an example of a numbered third-level heading.

%\subsubsection*{This is an unnumbered third-level section head}
%This is an example of an unnumbered third-level heading.

\section{Introduction}

It is well-known that the theory of chaos in finite-dimensional dynamical 
systems has been well-developed. That includes both discrete maps and systems 
of ordinary differential equations. Such theory has produced important 
mathematical theorems and led to important applications in physics, 
chemistry, biology, and engineering etc. \cite{GH83} \cite{LL83}. On the 
contrary, the theory of chaos in partial differential equations has not 
been well-developed. On the other hand, the demand for such a theory is 
much more stronger than for finite-dimensional systems. Mathematically, 
studies on infinite-dimensional systems pose much more challenging problems. 
For example, as phase spaces, Banach spaces possess much more structures 
than Euclidean spaces. In terms of applications, most of important natural
phenomena are described by partial differential equations, nonlinear wave 
equations, Yang-Mills equations, and Navier-Stokes equations, to name a few.

Nonlinear wave equations are the most important class of equations in 
natural sciences. They describe a wide spectrum of phenomena; motion of 
plasma, nonlinear optics (laser), water waves, vortex motion, to name a 
few. Among these nonlinear wave equations, there is a class of equations 
called soliton equations. This class of equations describes a variety of 
phenomena. In particular, the same soliton equation describes several 
different phenomena. For references, see for example \cite{AS81} 
\cite{AC91}. Mathematical theories on soliton equations have been well 
developed. Their Cauchy problems are completely solved through inverse 
scattering transforms. Soliton equations are integrable Hamiltonian 
partial differential equations which are the natural counterparts of 
finite-dimensional integrable Hamiltonian systems.

To set up a systematic study on chaos in partial differential equations, 
we started with the perturbed soliton equations. We 
classify the perturbed soliton equations into three categories: 
\begin{enumerate}
\item Perturbed (1+1)-Dimensional Soliton Equations, 
\item Perturbed Soliton Lattices, 
\item Perturbed ($1+n$)-Dimensional Soliton Equations ($n \geq 2$). 
\end{enumerate}
For each category, we chose a candidate to study. The integrable theories 
for every members in the same category are parallel, and for members in 
different categories are substantially different. The theorem on the 
existence of chaos for each candidate can be parallelly generalized to 
the rest members of the same category. 
\begin{itemize}
\item The candidate for Category 1 is the perturbed cubic focusing nonlinear 
Schr\"odinger equation \cite{LM94} \cite{LMSW96} \cite{Li99a},
\[
i\pa_t q = \pa^2_x q + 2 [ |q|^2 -\om^2 ]q + \mbox{ Perturbations},
\]
under periodic and even boundary conditions $q(x +1)=q(x)$ and $q(-x)=q(x)$,
$\om$ is a real constant.
\item The candidate for Category 2 is the perturbed discrete cubic focusing 
nonlinear Schr\"odinger equation \cite{Li92} \cite{LM97} \cite{LW97a}, 
\begin{eqnarray*}
i\dot{q_n}&=&{1 \over h^2}[q_{n+1}-2q_n+q_{n-1}]  \\
& & +|q_n|^2(q_{n+1}+q_{n-1})-2\om^2 q_n + \mbox{Perturbations} \ ,
\end{eqnarray*}
under periodic and even boundary conditions $q_{n+N}=q_n$ and $q_{-n}=q_n$.
\item The candidate for Category 3 is the perturbed Davey-Stewartson II 
equations \cite{Li00a}, 
\[
  \left \{ \begin{array}{l}
i \partial_t q = [ \partial^2_x - \partial^2_y]q+ [ 2(
|q|^2 - \omega^2) + u_y ] q + \mbox{Perturbations}, \cr \cr
[\partial^2_x + \partial^2_y] u = -4 \partial_y |q|^2 \, , \cr    
    \end{array} \right.
\]
under periodic and even boundary conditions
\begin{eqnarray*}
& &  q(t,x + l_x,y) = q(t,x,y) = q(t,x,y+l_y) \, , \\
& &  u(t,x + l_x,y) = u(t,x,y) = 
  u(t,x,y+l_y) \, ,
\end{eqnarray*}
and
\begin{eqnarray*}
& & q(t,-x,y) = q (t,x,y) = q(t,x,-y) \, ,\\
& & u(t,-x,y) = u(t,x,y) = u(t,x,-y) \, .
\end{eqnarray*}
\end{itemize}
We have established a standard program for 
proving the existence of chaos in perturbed soliton equations, with the 
machineries: 
\begin{enumerate}
\item Darboux Transformations for Soliton Equations.
\item Isospectral Theory for Soliton Equations Under Periodic Boundary 
Condition.
\item Persistence of Invariant Manifolds and Fenichel Fibers.
\item Melnikov Analysis.
\item Smale Horseshoes and Symbolic Dynamics Construction of Conley-Moser Type.
\end{enumerate}

The most important implication of the theory on chaos in partial 
differential equations in theoretical physics will be on the study of 
turbulence. For that goal, we chose the 2D Navier-Stokes 
equations under periodic boundary conditions to begin a dynamical system 
study on 2D turbulence. Since they possesses Lax pair \cite{Li01a}
and Darboux transformation \cite{LY00}, the 2D Euler equations are the 
starting point for an analytical study. The high Reynolds number 2D 
Navier-Stokes equations are viewed as a singular perturbation of the 2D Euler 
equations through the perturbation parameter $\e = 1/Re$ which is the 
inverse of the Reynolds number. Corresponding singular perturbations of 
nonlinear Schr\"odinger equation have been studied in \cite{Zen00a} 
\cite{Zen00b}
\cite{Li01b} \cite{Li01c}. We have studied the linearized 
2D Euler equations and obtained a complete spectra theorem \cite{Li00b}. 
In particular, we have identified unstable eigenvalues. Then we
found the approximate representations of the hyperbolic structures associated 
with the unstable eigenvalues through Galerkin truncations \cite{Li01d}. 

\section{Existence of Chaos in Perturbed Soliton Equations}

By existence of chaos, we mean that there exist a Smale horseshoe and 
the Bernoulli shift dynamics for certain Poincar\'e map. For lower dimensional 
systems, there have been a lot of theorems on the existence of chaos 
\cite{GH83} \cite{LL83}. For perturbed soliton equations under dissipative 
perturbations, we first establish the existence of a Silnikov 
homoclinic orbit. And then we define a Poincar\'e section which is 
transversal to the Silnikov homoclinic orbit, and the Poincar\'e map on the 
Poincar\'e section induced by the flow. Finally we construct the Smale 
horseshoe for the Poincar\'e map. In establishing the existence of the 
Silnikov homoclinic orbit, we need to build a Melnikov analysis 
through Darboux 
transformations to generate the explicit representation for the unperturbed 
heteroclinic orbit, the isospectral theory for soliton equations to generate 
the Melnikov vectors, and the persistence of invariant manifolds and Fenichel 
fibers. We also need to utilize the properties of the Fenichel fibers to 
build a second measurement inside a slow manifold, together with normal form 
techniques. The Melnikov measurement and the second measurement together 
lead to the existence of the Silnikov homoclinic orbit through implicit 
function arguments. In establishing the existence of Smale horseshoes for the 
Poincar\'e map, we first need to establish a smooth linearization in the 
neighborhood of the saddle point (i.e. the asymptotic point of the Silnikov 
homoclinic orbit). Then the dynamics in the neighborhood of the saddle 
point is governed by linear partial differential equations which are 
explicitly solvable. The global dynamics in the tubular neighborhood of 
the Silnikov homoclinic orbit away from the above neighborhood of the 
saddle point, can be approximated by linearized flow along the Silnikov 
homoclinic orbit due to finiteness of the passing time. Finally we can 
obtain a semi-explicit representation 
for the Poincar\'e map. Then we establish the existence of fixed points 
of the Poincar\'e map under certain except-one-point conditions. And we 
study the action of the Poincar\'e map in 
the neighborhood of these fixed points, and verify the Conley-Moser criteria 
to establish the existence of Smale horseshoes and Bernoulli shift dynamics. 

\subsection{Existence of Chaos in Perturbed (1+1)-Dimensional Soliton 
Equations}

For this category of the perturbed soliton equations, we chose the candidate 
to be the perturbed cubic nonlinear Schr\"odinger equation.
The cubic nonlinear Schr\"odinger equation describes self-focusing phenomena 
in nonlinear optics, deep water surface wave, vortex filament motion etc.. 
Recently, more and more interests are on perturbed nonlinear Schr\"odinger 
equations describing new nonlinear optical effects, for example, the works 
of the Laser Center at Oklahoma State University.

\subsubsection{Dissipative Perturbations}

In a series of three papers \cite{LM94} \cite{LMSW96} \cite{Li99a}, 
we proved the existence of chaos in the cubic nonlinear Schr\"odinger 
equation under dissipative perturbations. 
We study the following perturbed nonlinear Schr\"odinger equation:
\begin{equation}
iq_t = q_{xx} + 2 \bigg [ |q|^2-\om^2 \bigg ]q 
+ i\e \bigg [ -\al q +\hat{D}^2 q +\Ga \bigg ], \label{PNLS}
\end{equation}
under even periodic boundary conditions 
\[
q(-x)=q(x),\ \ \ \  q(x+1)=q(x); 
\]
where $i=\sqrt{-1}$, $q$ is a complex-valued function of two 
variables ($x,t$), ($\om,\al,\Ga$) are positive constants, $\e$ is the 
positive perturbation parameter, $\hat{D}^2$ is a ``regularized'' Laplacian 
specifically defined by
\[
\hat{D}^2 q \equiv -\sum^{\infty}_{j=1} \be_j k_j^2 \hat{q}_j \cos k_j x,
\]
in which $k_j=2\pi j$, $\hat{q}_j$ is the Fourier transform of $q$, 
$\be_j=\be$ for $j \leq N$, $\be_j = \al_* k_j^{-2}$ for $j > N$, 
$\be$ and $\al_*$ are positive constants, and $N$ is a large fixed 
positive integer.
\begin{theorem}[Homoclinic Orbit Theorem]
There exists a positive number $\e_0$ such that for any $\e \in (0,\e_0)$,
there exists a codimension 1 hypersurface $E_\e$ in the external parameter
space $\{ \om, \al, \Ga, \be, \al_* \}$. For any external parameters 
$(\om, \al, \Ga, \be, \al_*) \in E_\e$, there exists a symmetric pair of 
homoclinic orbits $h_k =h_k(t,x)\ \ (k=1,2)$ in $H^1_{e,p}$ 
(the Sobolev space $H^1$ of even and periodic functions) for the PDE
(\ref{PNLS}), which are asymptotic to a fixed point $q_\e$. The symmetry
between $h_1$ and $h_2$ is reflected by the relation that $h_2$ is a 
half-period translate of $h_1$, i.e. $h_2(t,x)=h_1(t,x+1/2)$. The 
hypersurface $E_\e$ is a perturbation of a known surface $\be =\k(\om) \al$,
where $\k(\om)$ is shown in Figure \ref{pdek}.
\label{homthm}
\end{theorem}
\begin{figure}
\vspace{4.0in}
\caption{The curve of $\k = \k(\om)$.}
\label{pdek}
\end{figure}
For the complete proof of the theorem, see \cite{LM94} and \cite{LMSW96}.
The main argument is a combination of a Melnikov analysis and a geometric 
singualr perturbation theory for partial differential equations. The Melnikov
function is evaluated on a homoclinic orbit of the nonlinear Schr\"odinger 
equation, generated through Darboux transformations. For more details on 
this, see the later section on the Darboux transformations for the discrete 
nonlinear Schr\"odinger equation.
\begin{theorem}[Horseshoe Theorem]
Under certain generic assumptions for the perturbed nonlinear
Schr{\"{o}}dinger system (\ref{PNLS}), there exists a compact Cantor subset 
$\La$ of $\Sg$ (a Poincar\'e section transversal to the homoclinic orbit), 
$\La$ consists of points, and is invariant under $P$ (the Poincar\'e map 
induced by the flow on $\Sg$).
$P$ restricted to $\La$, is topologically conjugate to the shift 
automorphism $\chi$ on four symbols $1, 2, -1, -2$. That is, there exists
a homeomorphism
\[
\phi \ : \ \W \mapsto  \La,
\]
\nid
where $\W$ is the topological space of the four symbols, such that the 
following diagram commutes:
\begin{equation} 
\begin{array}{ccc}
\W &\maprightu{\phi} & \Lambda \\
\mapdownl{\chi} & & \mapdownr{P} \\
\W & \maprightd{\phi} & \Lambda
\end{array} 
\nonumber
\end{equation}
\label{horseshthm}
\end{theorem}
For the complete proof of the theorem, see \cite{Li99a}.
The construction of horseshoes is of Conley-Moser type for partial 
differential equations.

\subsubsection{Singular Perturbations}

Recently, singular perturbation, i.e. replacing $\hat{D}^2 q$ by $\pa^2_x q$,
has been studied \cite{Zen00a} \cite{Zen00b} \cite{Li01b} \cite{Li01c}. 
Consider the singularly perturbed nonlinear Schr\"odinger equation,
\begin{equation}
iq_t = q_{xx} +2 [ |q|^2 - \om^2] q +i \e [- \al q + \be q_{xx}+ \Ga ] \ ,
\label{spnls}
\end{equation}
where $q = q(t,x)$ is a complex-valued function of the two real 
variables $t$ and $x$, $t$ represents time, and $x$ represents
space. $q(t,x)$ is subject to periodic boundary condition of period 
$1$, and even constraint, i.e. 
\[
q(t,x + 1) = q(t,x)\ , \ \ q(t,-x) = q(t,x)\ .
\]
$\om$ is a positive constant, $\al >0$, $\be >0$, and $\Ga$ are constants, 
and $\e > 0$ is the perturbation parameter.
The main difficulty introduced by the singular 
perturbation $\e \pa_x^2$ is that it breaks the spectral gap condition 
of the unperturbed system. Therefore, standard invariant manifold results 
will not apply. Nevertheless, it turns out that certain invariant manifold 
results do hold. The regularity of such invariant manifolds at $\e = 0$ is 
controled by the regularity of $e^{\e \pa_x^2}$ at $\e = 0$. 
\begin{theorem}[Li, \cite{Li01b}] 
There exists a $\e_0 > 0$, such that for any $\e \in (0, \e_0)$, there exists 
a codimension 1 surface $E_\e$ in the space of $(\om, \al, \be, \Ga) \in 
R^+\times R^+\times R^+\times R^+$,
where $\om \in (\pi, 2\pi)/S$, $S$ is a finite subset.
For any external parameters on the codimension-one
surface, the perturbed nonlinear Schr\"odinger equation 
(\ref{spnls}) possesses a symmetric pair of 
homoclinic orbits $h_k =h_k(t,x)\ \ (k=1,2)$ in $C^\infty_{e,p}[0,1]$
(the space of $C^\infty$ even and periodic functions on the interval [0,1]) 
, which is asymptotic to a saddle fixed point $q_\e$. The symmetry
between $h_1$ and $h_2$ is reflected by the relation that $h_2$ is a 
half-period translate of $h_1$, i.e. $h_2(t,x)=h_1(t,x+1/2)$. The 
hypersurface $E_\e$ is a perturbation of a known surface $\be =\k(\om) \al$,
where $\k(\om)$ is shown in Figure \ref{pdek}.
\end{theorem}

\subsubsection{Hamiltonian Perturbations}

The problem on the existence of chaos in the cubic nonlinear Schr\"odinger 
equations under Hamiltonian perturbations is largely open. The right 
objects to investigate should be ``homoclinic tubes'' rather than 
``homoclinic orbits'' due to the non-dissipative nature and 
infinite-dimensionality of the perturbed system. Transversal 
homoclinic tubes are objects of large dimensional generalization of 
transversal homoclinic orbits. As Smale's theorem indicates, structures 
in the neighborhood of a transversal homoclinic orbit are rich, structures 
in the neighborhood of a transversal homoclinic tube are even richer. 
Especially in high dimensions, dynamics inside each invariant 
tubes in the neighborhoods of homoclinic tubes are often 
chaotic too. We call such chaotic dynamics ``{\em{chaos in 
the small}}'', and the symbolic dynamics of the invariant 
tubes ``{\em{chaos in the large}}''. Such cascade structures 
are more important than the structures in a neighborhood of a 
homoclinic orbit, when high or infinite dimensional dynamical 
systems are studied. Symbolic dynamics structures in the 
neighborhoods of homoclinic tubes are more observable than 
in the neighborhoods of homoclinic orbits in numerical and 
physical experiments. When studying high or infinite 
dimensional Hamiltonian system (for example, the cubic nonlinear 
Schr\"odinger equation under Hamiltonian perturbations), each 
invariant tube contains both KAM tori and stochastic layers 
(chaos in the small). Thus, not only dynamics inside each 
stochastic layer is chaotic, all these stochastic layers also 
move chaotically under Poincar\'e maps.

There have been a lot 
of works on the KAM theory of soliton equations under Hamiltonian 
perturbations \cite{Way90} \cite{CW93} \cite{Kuk93} \cite{Bou94} 
\cite{Pos96}. For perturbed nonlinear Schr\"odinger equations, 
we are interested in the region of the phase space where there exist 
hyperbolic structures. Thus, the relevant KAM tori are hyperbolic. In 
finite dimensions, the relevant work on such tori is that of Graff 
\cite{Gra74}. In infinite dimensions, the author is not aware of such work
yet.

In the paper \cite{Li99b}, the author studied the cubic nonlinear 
Schr\"odinger equation under Hamiltonian perturbations:
\begin{equation}
i q_t = q_{xx} + 2 [|q|^2-\om^2] q + \e [\al_1 +2 \al_2 \bq ]\ ,
\label{hpnls}
\end{equation}
under even periodic boundary conditions $q(-x)=q(x)$ and $q(x+1)=q(x)$;
where $i = \sqrt{-1}$, $q$ is a complex-valued function of two variables 
($t,x$), ($\om,\al_1,\al_2$) are real constants, $\e$ is the perturbation 
parameter. The system (\ref{hpnls}) can be written in the Hamiltonian form:
\[
iq_t =\dl H / \dl \bq , 
\]
where $H= H_0 + \e H_1 $,
\[
H_0 = \int^1_0 [ |q|^4 - 2 \om^2 |q|^2 -|q_x|^2] dx ,
\]
\[
H_1 = \int^1_0 [ \al_1 (q+ \bq ) +\al_2 (q^2 +\bq^2 )] dx .
\]
\begin{definition}
Denote by $W^{(c)}$ a normally hyperbolic 
center manifold, by $W^{(cu)}$ and $W^{(cs)}$ the center-unstable 
and center-stable manifolds such that $W^{(c)}=W^{(cu)} \cap 
W^{(cs)}$, and by $F^t$ the evolution operator of the partial 
differential equation. Let $\HH$ be a submanifold in the 
intersection between the center-unstable and center-stable 
manifolds $W^{(cu)}$ and $W^{(cs)}$, such that for any point 
$q \in \HH$, distance$\{ F^t(q),W^{(c)}\} \ra 0$, as $|t| \ra 
\infty$. We call $\HH$ a transversal homoclinic tube asymptotic 
to $W^{(c)}$ under the flow $F^t$ if the intersection between 
$W^{(cu)}$ and $W^{(cs)}$ is transversal at $\HH$. Let $\Sg$ be 
an Poincar\'e section which intersects $\HH$ transversally, and $P$ is 
the Poincar\'e map induced by the flow $F^t$; then $\HH \cap \Sg$ is called a 
transversal homoclinic tube under the Poincar\'e map $P$.
\end{definition}
\begin{theorem}[Homoclinic Tube Theorem]
There exist a positive constant $\e_0 > 0$ and a region $\E$ for ($\al_1,
\al_2,\om$), 
such that for any $\e \in (-\e_0,\e_0)$ and any 
$(\al_1,\al_2,\om) \in \E$, there 
exists a codimension 2 transversal homoclinic tube asymptotic to 
a codimension 2 center manifold $W^{(c)}$.
\label{thmht}
\end{theorem}
For a complete proof of this theorem, see \cite{Li99b}.

\subsection{Chaos in Perturbed Soliton Lattices} 

For this category, we chose the candidate to be the perturbed cubic 
nonlinear Schr\"odinger lattice.

\subsubsection{Dissipative Perturbations} 

In a series of three papers \cite{Li92} \cite{LM97} \cite{LW97a}, we 
proved the existence of chaos in the discrete cubic nonlinear Schr\"odinger 
equation under a concrete dissipative perturbation. 

We study the perturbed discrete cubic nonlinear Schr\"odinger equation
\begin{eqnarray}
i\dot{q}_n &=& {1\over h^2} \bigg [ q_{n+1}-2 q_n +q_{n-1} \bigg ]
+|q_n|^2(q_{n+1}+q_{n-1})-2\om^2 q_n \nonumber \\
&+& i\e \bigg [ -\al q_n +{\be \over h^2}(q_{n+1}-2 q_n +q_{n-1})+\Ga 
\bigg ], \label{PDNLS}
\end{eqnarray}
under even periodic boundary conditions ($q_{N-n}=q_n$) and ($q_{n+N}=q_n$)
for arbitrary $N$; where $i=\sqrt{-1}$, $q_n$'s are complex variables, 
$h=1/N$, ($\om,\al,\be,\Ga$) are positive constants, $\e$ is the positive 
perturbation parameter.

Denote by $\Sg_N\ (N\geq 7)$ the external parameter space,
\begin{eqnarray*}
\Sg_N&=&\bigg\{ (\om,\al,\be,\Ga)\ \bigg | \ \om \in (N\tan{\pi \over N}, 
               N\tan{2\pi \over N}),\\ 
     & &\Ga\in (0,1), \al\in (0,\al_0), \be\in (0,\be_0); \\
     & &\mbox{where}\ \al_0\ \mbox{and}\ \be_0\ \mbox{are}\ \mbox{any}
        \ \mbox{fixed}\ \mbox{positive}\ \mbox{numbers}. \bigg\}
\end{eqnarray*}
\begin{theorem}
For any $N$ ($7\leq N<\infty$), there exists a positive number $\e_0$, 
such that for any $\e \in (0,\e_0)$, there exists a 
codimension $1$ submanifold $E_\e$ in $\Sg_N$; for any
external parameters ($\om,\al,\be,\Ga$) on $E_\e$, there exists a 
homoclinic orbit asymptotic to a fixed point $q_\e$.
The submanifold $E_\e$ is 
in an $O(\e^\nu)$ 
neighborhood of the hyperplane $\be=\k\ \al$, where $\k=\k(\om;N)$
is shown in Figures \ref{disk1} and \ref{disk2}, 
$\nu=1/2 -\dl_0, \ 0<\dl_0 \ll 1/2$. 
\end{theorem}
\begin{figure}
\vspace{7.5in}
\caption{The curve of $\k = \k(\om;N)$ .}
\label{disk1}
\end{figure}
\begin{figure}
\vspace{7.5in}
\caption{The curve of $\k = \k(\om;N)$ .}
\label{disk2}
\end{figure}
\begin{remark}
In the cases ($3\leq N \leq 6$), $\k$ is always negative as shown in 
Figure \ref{disk2}. Since we require both dissipation parameters 
$\al$ and $\be$ 
to be positive, the relation $\be=\k \al$ shows that the existence 
of homoclinic orbits violates this positivity. For $N \geq 7$, $\k$ 
can be positive as shown in Figure \ref{disk1}. When $N$ is even and 
$\geq 7$, there is in fact a pair of homoclinic orbits asymptotic to 
a fixed point $q_\e$ at the same values of the external parameters; 
since for even $N$, we have the symmetry:
If $q_n=f(n,t)$ solves (\ref{PDNLS}),
then $q_n=f(n+N/2,t)$ also solves (\ref{PDNLS}). When $N$ is odd
and $\geq 7$, the study can not guarantee that two homoclinic orbits 
exist at the same value of the external parameters.
\end{remark}
For the complete proof of this theorem, see \cite{LM97}.
\begin{theorem}[Horseshoe Theorem]
Under certain generic assumptions for the perturbed discrete nonlinear
Schr{\"{o}}dinger system (\ref{PDNLS}), there exists a compact Cantor subset 
$\La$ of $\Sg$ (a Poincar\'e section transversal to the homoclinic orbit), 
$\La$ consists of points, and is invariant under $P$ (the Poincar\'e map 
induced by the flow on $\Sg$).
$P$ restricted to $\La$, is topologically conjugate to the shift 
automorphism $\chi$ on four symbols $1, 2, -1, -2$. That is, there exists
a homeomorphism
\[
\phi \ : \ \W \mapsto  \La,
\]
\nid
where $\W$ is the topological space of the four symbols, such that the 
following diagram commutes:
\begin{equation} 
\begin{array}{ccc}
\W &\maprightu{\phi} & \Lambda \\
\mapdownl{\chi} & & \mapdownr{P} \\
\W & \maprightd{\phi} & \Lambda
\end{array} 
\nonumber
\end{equation}
\end{theorem}
For the complete proof of the theorem, see \cite{LW97a}.

The unperturbed homoclinic orbits for the discrete nonlinear Schr\"odinger 
equation
\begin{equation}
i\dot{q_n}={1 \over h^2}\bigg[q_{n+1}-2q_n+q_{n-1}\bigg]+|q_n|^2(q_{n+1}+
           q_{n-1})-2\om^2 q_n, \label{IDNLS}
\end{equation}
was constructed through the Darboux transformations which will be 
presented below in details.
The discrete nonlinear Schr\"odinger equation is associated with the 
following discrete Zakharov-Shabat system \cite{AL76}:
\begin{eqnarray}
\varphi_{n+1}&=&L^{(z)}_n\varphi_n, \label{ZS1} \\
\dot{\varphi}_n&=&B^{(z)}_n\varphi_n, \label{ZS2}
\end{eqnarray}
\noindent
where
\begin{eqnarray*}
L^{(z)}_n&\equiv&\left( \begin{array}{cc} 
                                            z& ihq_n \cr
                                            ih\bq_n & 1/z \cr
                                       \end{array} \right), \\
B^{(z)}_n&\equiv&{i\over h^2}\left( \begin{array}{cc}
  1-z^2+2i\la h-h^2q_n\bq_{n-1}+\om^2 h^2 & -zihq_n+(1/z)ihq_{n-1}  \\
  -izh\bq_{n-1}+(1/z)ih\bq_n& 1/z^2-1+2i\la h+h^2\bq_nq_{n-1}-\om^2 h^2
                                 \end{array} \right),
\end{eqnarray*}
and where $z\equiv \exp(i\la h)$.

Fix a solution $q_n(t)$ of the 
system (\ref{IDNLS}), for which the linear operator $L_n$ has a double 
point $z^d$ of geometric multiplicity 2, which is not on the unit circle.
We denote two linearly independent solutions (Bloch functions) of the 
discrete Lax pair (\ref{ZS1};\ref{ZS2}) at $z=z^d$ by $(\phi_n^+,\phi_n^-)$.
Thus, a general solution of the discrete Lax pair (\ref{ZS1};\ref{ZS2}) 
at $(q_n(t),z^d)$ is given by 
\[
\phi_n(t; z^d, c)=\phi_n^+ + c \phi_n^-,
\]
\nid
where $c$ is a complex parameter called B\"acklund parameter. We use $\phi_n$ 
to define a transformation matrix $\Ga_n$ by
\[
\Ga_n=\left(\begin{array}{cc} z+(1/z)a_n & b_n \cr c_n &-1/z+z d_n \cr
            \end{array} \right),
\]
\nid
where,
\begin{eqnarray*}
a_n &=& {z^d \over (\bar{z}^d)^2\Dl_n}\bigg [|\phi_{n2}|^2+|z^d|^2|\phi_{n1}|^2
    \bigg ],\\
d_n &=& -{1 \over z^d\Dl_n}\bigg [|\phi_{n2}|^2+|z^d|^2|\phi_{n1}|^2
    \bigg ],\\
b_n &=& {|z^d|^4-1 \over (\bar{z}^d)^2\Dl_n}\phi_{n1}\bar{\phi}_{n2}, \\
c_n &=& {|z^d|^4-1 \over z^d\bar{z}^d\Dl_n}\bar{\phi}_{n1}\phi_{n2}, \\
\Dl_n &=& -{1 \over \bar{z}^d}\bigg [|\phi_{n1}|^2+|z^d|^2|\phi_{n2}|^2
    \bigg ].
\end{eqnarray*}
Then we define $Q_n$ and $\Psi_n$ by
\begin{equation}
Q_n\equiv {i\over h}b_{n+1}-a_{n+1}q_n
\label{BD1}
\end{equation}
\nid
and
\begin{equation}
\Psi_n(t;z)\equiv \Ga_n(z;z^d;\phi_n)\psi_n(t;z)
\label{BD2}
\end{equation}
\nid
where $\psi_n$ solves the discrete Lax pair (\ref{ZS1};\ref{ZS2}) 
at $(q_n(t),z)$. Formulas (\ref{BD1}) and (\ref{BD2}) are the 
B\"acklund-Darboux transformations for the potential and eigenfunctions, 
respectively. We have the following theorem \cite{Li92}.
\begin{theorem}[B\"acklund-Darboux Transformations]
Let $q_n(t)$ denote a solution of the system (\ref{IDNLS}), for which 
the linear operator $L_n$ has a double point $z^d$ of geometric 
multiplicity 2, 
which is not on the unit circle and which is associated with an 
instability. We denote two linearly independent 
solutions (Bloch functions) of the discrete Lax pair (\ref{ZS1};\ref{ZS2}) 
at $(q_n, z^d)$ by $(\phi_n^+,\phi_n^-)$. We define $Q_n(t)$ and $\Psi_n(t;z)$ 
by (\ref{BD1}) and (\ref{BD2}). Then
\begin{enumerate}
\item $Q_n(t)$ is also a solution of the system (\ref{IDNLS}). (The eveness
of $Q_n$ can be guaranteed by choosing the complex B\"acklund parameter $c$ to 
lie on an certain curve.)
\item $\Psi_n(t;z)$ solves the discrete Lax pair (\ref{ZS1};\ref{ZS2}) 
at $(Q_n(t),z)$.
\item $\Dl(z;Q_n)=\Dl(z;q_n)$, for all $z\in C$, where $\Dl$ is the 
Floquet discriminant.
\item $Q_n(t)$ is homoclinic to $q_n(t)$ in the sense that $Q_n(t) \ra 
e^{i\th_{\pm}}\ q_n(t)$, exponentially as $\exp (-\sg |t|)$ as $t \ra 
\pm \infty$. Here $\th_{\pm}$ are the phase shifts, $\sg$ is a 
nonvanishing growth rate associated to the double point $z^d$, and 
explicit formulas can be developed for this growth rate and 
for the phase shifts $\th_{\pm}$.
\end{enumerate}
\label{Backlund}
\end{theorem}
Next we consider a concrete example. Let
\begin{equation}
q_n=q,\ \forall n; \ \ q=a\exp\{-2i[(a^2-\om^2)t]+i\ga \}, 
\label{ucsl}
\end{equation}
where $N\tan{\pi \over N}< a <
N\tan{2\pi \over N}\ \mbox{for}\ N>3, 3\tan{\pi \over 3}< a < 
\infty\ \mbox{for}\ N=3 $. 
Then $Q_n$ defined in (\ref{BD1})
has the explicit representation:
\begin{equation}
Q_n\equiv Q_n(t;\ N,\om,\ga,r,\pm)=q\bigg [{G\over H_n}-1\bigg ],
\label{horu}
\end{equation}
\nid
where,
\[
G=1+\cos 2P-i\sin 2P \tanh\tau,
\]
\[
H_n=1\pm{1 \over \cos \vth}\sin P \ \mbox{sech}\ \tau \cos 2n\vth,
\]
\[
\tau=4N^2\sqrt{\rho}\sin \vth\sqrt{\rho \cos^2\vth-1}\ t+r,
\]
\nid
where $r$ is a real parameter. Furthermore,
\[
P=\arctan{{\sqrt{\rho \cos^2\vth-1}\over \sqrt{\rho}\sin \vth}},
\]
\[
\vth={\pi \over N}, \ \rho=1+{|q|^2 \over N^2}.
\]
\nid
As $\tau\ra \pm \infty$, $Q_n\ra qe^{\mp i2P}$. Therefore, $Q_n$ is 
homoclinic to the circle $|q_n|=a$, and heteroclinic to points on 
the circle which are separated in phase of $-4P$.

\subsubsection{Hamiltonian Perturbations} 

In the paper \cite{Li98a}, the author studied the discrete 
nonlinear Schr\"odinger equation under Hamiltonian perturbations:
\begin{eqnarray}
i \dot{q}_n &=& {1 \over h^2} [q_{n+1} -2 q_n +q_{n-1}] + |q_n|^2(q_{n+1} 
+q_{n-1}) - 2 \om^2 q_n 
\label{hpdnls} \\
&+& \e \bigg \{ [\al_1 (q_n +\bq_n) +\al_2 (q_n^2 +\bq_n^2)]q_n +
[\al_1 +2 \al_2 \bq_n] {\rho_n \over h^2} \ln \rho_n \bigg \}\ , \nonumber 
\end{eqnarray}
where $i = \sqrt{-1}$, $q_n's$ are complex variables, $n \in Z$, 
($\om, \al_1, \al_2$) are real constants, $\e$ is the perturbation 
parameter, $h$ is the step size, $h =1/N$, $N \geq 3$ is an integer,
$\rho_n =1 +h^2 |q_n|^2$, and $q_{n+N}=q_n\ , \ \ \ q_{-n}=q_n$.
The system (\ref{hpdnls}) can be written in the Hamiltonian form:
\[
i \dot{q}_n = \rho_n \ {\pa H \over \pa \bq_n} ,
\]
where $H=H_0 + \e H_1 $,
\[
H_0 = {1 \over h^2} \sum_{n=0}^{N-1} [ \bq_n (q_{n+1} +q_{n-1}) - 
{2 \over h^2} (1+\om^2 h^2) \ln \rho_n ] , 
\]
\[
H_1 = {1 \over h^2} \sum_{n=0}^{N-1} [ \al_1 (q_n + \bq_n) +\al_2 
(q_n^2 +\bq_n^2)] \ln \rho_n .
\]
\begin{theorem}[Homoclinic Tube Theorem]
There exist a positive constant $\e_0 > 0$ and a region $\E$ for 
($\al_1,\al_2,\om$), 
such that for any $\e \in (-\e_0,\e_0)$ and any $(\al_1,\al_2,\om) 
\in \E$, there 
exists a codimension 2 transversal homoclinic tube asymptotic to 
a codimension 2 center manifold $W^{(c)}$.
\end{theorem}
For a complete proof of this theorem, see \cite{Li98a}.

\subsection{Chaos in Perturbed ($1+n$)-Dimensional Soliton 
Equations ($n \geq 2$)}

For this category of the perturbed soliton equations, we chose the 
candidate to be the perturbed Davey-Stewartson II equations.
The Davey-Stewartson II equations describe nearly one-dimensional water 
surface wave train \cite{DS74}. There have been a lot of studies on the 
inverse scattering transforms for this set of equations \cite{AC91} 
\cite{AS81}. The inverse scattering transforms for (1+n)-dimensional 
soliton equations ($n \geq 2$) are substantially different from those 
for (1+1)-dimensional soliton equations and soliton lattices. In fact, 
the Davey-Stewartson II equations possess finite-time singularities 
\cite{Oza92}. For the perturbed Davey-Stewartson II equations, the theory 
on chaos is largely unfinished. So far, its Melnikov theory has been 
successfully built.

Although the inverse spectral theory for the DSII equations is very 
different from those for (1+1)-dimensional soliton equations and there 
is no Floquet spectral theory, its B\"acklund-Darboux transformation is 
as simple as those for (1+1)-dimensional soliton equations, e.g. the 
cubic nonlinear Schr\"odinger equation. 
These B{\"{a}}cklund-Darboux transformations are
successfully utilized to construct heteroclinic orbits of
Davey-Stewartson II equations through an elegant iteration of the
transformations. In \cite{LM94}, we successfully built Melnikov
vectors for the focusing cubic nonlinear Schr\"{o}dinger equation
with the gradients of the invariants $F_j$ defined through the
Floquet discriminants evaluated at critical spectral points. The invariants
$F_j$'s Poisson commute with the Hamiltonian, and their gradients
decay exponentially as time approaches positive and negative infinities --
these two properties are crucial in deriving and evaluating Melnikov functions.
Since there
is no Floquet discriminant for Davey-Stewartson equations (in
contrast to nonlinear Schr{\"{o}}dinger equations \cite{LM94}),
the Melnikov vectors here are built with the novel idea of replacing
the gradients of Floquet discriminants by quadratic products of Bloch
functions. Such Melnikov vectors still maintain the properties of
Poisson commuting with the gradient of the Hamiltonian and exponential
decay as time approaches positive and negative infinities. This solves
the problem of building Melnikov vectors for Davey-Stewartson equations
without using the gradients of Floquet discriminant. Melnikov functions
for perturbed Davey-Stewartson II equations evaluated on the above
heteroclinic orbits are built. 

\subsubsection{Darboux Transformations}

First we study the Darboux transformations for the Davey-Stewartson II 
(DSII) equations:
\begin{eqnarray}
& & i \pa_t q = [\pa_x^2 - \pa_y^2] q + [2(|q|^2-\om^2) + u_y] q\ , 
\label{DS2} \\ 
& & \ \ [\pa_x^2 + \pa_y^2] u = - 4 \pa_y |q|^2 \ ; \nonumber 
\end{eqnarray}
under periodic boundary conditions $q(t,x+l_x,y)=q(t,x,y+l_y)=q(t,x,y)$,
where $q$ and $u$ are a complex-valued 
and a real-valued functions of three variables ($t,x,y$). To simplify the 
study, we may also pose even conditions in both $x$ and $y$. The DSII 
equations are associated with a Lax pair and a congruent Lax pair.
The Lax pair is:
\begin{eqnarray}
         L \psi &=& \lambda \psi\,, \label{LP1} \\
       \partial_t \psi &=& A \psi\,,\label{LP2}
\end{eqnarray}
where $\psi = \left( \psi_1, \psi_2\right)^T$, and 
\[
  L = \left(
\begin{array}{lr}
D^{-} & q\\ \\
r & D^{+}
\end{array}
\right)\,,
\]
\[
A = i \left[
2 \left(
\begin{array}{cc}
- \partial^2_x & q \partial_x\\
r \partial_x & \partial^2_x
\end{array}
\right) \, + \,
\left(
\begin{array}{cc}
r_1 & (D^+ q)\\
-(D^{-} r) & r_2
\end{array}
\right)
\right]\,.
\]
Here we denote by
\begin{equation}
  D^+ = \alpha \partial_y + \partial_x\,, \qquad D^{-} = \alpha
  \partial_y - \partial_x\,.
\label{DD}
\end{equation}
where $r = \bar{q}$, $\al^2 = -1$, 
\[
  r_1 = \frac{1}{2} [-U+iV] \ , \ \ \  r_2= \frac{1}{2} [U+iV] \, , \ \ \ 
U=2(|q|^2 - \om^2 )+u_y.
\]
The congruent Lax pair is:
\begin{eqnarray}
 \hat{L} \hat{\psi} &=& \lambda \hat{\psi}
                               \,,\label{CLP1} \\
 \partial_t \hat{\psi} &=& \hat{A}
                                   \hat{\psi}
                                 \,,\label{CLP2}
\end{eqnarray}
where $\hat{\psi} = (\hat{\psi}_1, \hat{\psi}_2)^T$, and
\[
  \hat{L} =
  \left(
    \begin{array}{cc}
- D^+ & q\\ \\
r & -D^-
    \end{array}
  \right)\,,
\]
\[
\hat{A} = i \left[
  2 \left(
    \begin{array}{cc}
- \partial^2_x & q \partial_x\\
r \partial_x & \partial^2_x
    \end{array}
  \right) + 
  \left(
    \begin{array}{cc}
-r_2 & -(D^- q)\\
(D^+r) & -r_1
    \end{array}
  \right)
\right]\,.
\]
Let $(q,r= \bar{q}, r_1, r_2)$ be a solution to the DSII equation,
and let $\lambda_0$ be any value of $\lambda$. Denote by
$\psi = (\psi_1, \psi_2)^T$ the eigenfunction solving the Lax
pair (\ref{LP1}, \ref{LP2}) at $(q, r = \bar q, r_1, r_2;
\lambda_0)$. Define the matrix operator:
\begin{displaymath}
  \Gamma = 
\left[
  \begin{array}{cc}
             \wedge + a & b\\
             c & \wedge + d
  \end{array}
\right]\,,
\end{displaymath}
where $\wedge = \alpha \partial_y - \lambda$, and $a$, $b$, $c$,
$d$ are functions defined as:
\begin{eqnarray*}
  a &=& \frac{1}{\Delta} \left[ \psi_2 \wedge_2 \bar{\psi}_2 +
                  \beta \bar{\psi}_1 \wedge_1 \psi_1 \right]\,,\\[2ex]
  b &=& \frac{1}{\Delta} \left[ \bar{\psi}_2 \wedge_1 \psi_1 -
                   \psi_1 \wedge_2 \bar{\psi}_2 \right]\,,\\[2ex]
  c &=& \frac{\beta}{\Delta} \left[ \bar{\psi}_1 \wedge_1 \psi_2
                     - \psi_2 \wedge_2 \bar{\psi}_1 \right]\,,\\[2ex]
  d &=&  \frac{1}{\Delta} \left[ \bar{\psi}_2 \wedge_1 \psi_2 +
                      \beta \psi_1 \wedge_2 \bar{\psi}_1 \right]\,,
\end{eqnarray*}
in which $\wedge_1 = \alpha \partial_y - \lambda_0$, $\wedge_2 =
\alpha \partial_y + \bar{\lambda}_0$, and
\begin{displaymath}
  \Delta = - \left[ \beta | \psi_1 |^2 + |\psi_2|^2 \right]\,.
\end{displaymath}
Define a transformation as follows:
\begin{displaymath}
  \left\{
    \begin{array}{ccc}
(q,r=\beta \bar{q},r_1,r_2) &\rightarrow& (Q,R,R_1,R_2)\,, \\
\phi &\rightarrow& \Phi\,;
    \end{array}
\right.
\end{displaymath}
\begin{eqnarray}
                    Q   &=& q - 2b\,,\nonumber \\[2ex]
                    R   &=& \beta \bar{q} - 2c\,,\nonumber\\[2ex]
                    R_1 &=& r_1 + 2(D^+a)\,, \label{DSBT}\\[2ex]
                    R_2 &=& r_2 - 2 (D^- d)\,,\nonumber\\[2ex]
                   \Phi &=& \Gamma \phi\,;\nonumber
\end{eqnarray}
where $\phi$ is an eigenfunction solving the Lax pair (\ref{LP1},
\ref{LP2}) at $(q, r = \bar{q}, r_1, r_2; \lambda)$, $D^+$
and $D^-$ are defined in (\ref{DD}), 
\begin{theorem}[\cite{Li00a}]
The transformation (\ref{DSBT}) is a B\"acklund-Darboux
transformation. That is, the functions $(Q, R=\bar{Q}, R_1, R_2)$ defined
through the transformation (\ref{DSBT}) are also a solution to the
Davey-Stewartson II equations.  
The function $\Phi$ defined through the transformation
(\ref{DSBT}) solves the Lax pair (\ref{LP1}, \ref{LP2}) at $(Q,
R=\bar{Q}, R_1, R_2; \lambda)$.
\label{DSTH}
\end{theorem}
A concrete example with two iterations of the Darboux transformations
has been worked out in \cite{Li00a}.

\subsubsection{Melnikov Vectors}

The DSII equations can be put into the Hamiltonian form,

\begin{equation}
  \left\{
  \begin{array}{ccc}
    i q_t &=& \delta H / \delta \overline{q}  \ ,\cr
   i \overline{q}_t &=& - \delta H / \delta q \ , \cr 
  \end{array} \right. 
\label{fhDS2}
\end{equation}
where
\begin{displaymath}
  H= \int^{l_y}_0 \int^{l_x}_0 
  [\left| q_y \right|^2 - \left| q_x \right|^2 + 
       \frac{1}{2} (r_2-r_1) \left| q \right|^2] \, dx \, dy \, .
\end{displaymath}
Let $\psi = (\psi_1 , \psi_2)^T$ be an eigenfunction solving the
Lax pair (\ref{LP1}, \ref{LP2}), and $\hat{\psi} = (\hat{\psi}_1 ,
\hat{\psi}_2)^T$ be an eigenfunction solving the corresponding 
congruent Lax pair (\ref{CLP1}, \ref{CLP2}); then
\begin{lemma}
The inner product of the vector
\begin{equation}
  \U= \left(
    \begin{array}{c}
      \psi_2  \hat{\psi}_2 \cr
      \psi_1  \hat{\psi}_1
    \end{array}\right)^- +S \left(
    \begin{array}{c}
      \psi_2  \hat{\psi}_2 \\
      \psi_1  \hat{\psi}_1
    \end{array}\right) \, ,
\nonumber
\end{equation}
where $S= \displaystyle{\left(
    \begin{array}{ccc}
      0 & 1 \\ 1 & 0
    \end{array}
\right)}$, with the vector field $J \na H$ given by the 
right hand side of (\ref{fhDS2}) vanishes,
\begin{displaymath}
  \langle \U\, , \,  J \na H \rangle =0 \, .
\end{displaymath}
where
\begin{displaymath}
  \langle f,g \rangle = \int^{l_y}_0 \int^{l_x}_0
  \left\{ \overline{f}_1 g_1 + \overline{f}_2 g_2 \right\}
  \, dx \, dy \, .
\end{displaymath}
and 
\[
J = \left(
  \begin{array}{ll}
0 & 1\\
1 & 0
  \end{array}
\right).
\]
\label{melem2}
\end{lemma}

Consider the perturbed DSII equations
\begin{equation}
\left\{
    \begin{array}{l}
      i \partial_t q  =  [\partial^2_x - \partial^2_y] q + 
      [2 ( \left| q \right|^2 - \om^2) + u_y] q + 
      \e i f \, , \\[1ex]
[\partial^2_x + \partial^2_y] u 
       = -4 \partial_y \left| q \right|^2 \, ,
    \end{array} \right.
\label{PDS2}
\end{equation}
where $f$ is the perturbation which can
depend on $q$ and $\overline{q}$ and their derivatives and $t$, $x$
and $y$. Let $\vec{G}= (f, \overline{f})^T$. Then the Melnikov function 
has the expression,
\begin{eqnarray}
  M &=&  \int^{\infty}_{- \infty} \langle \U ,
                 \vec{G} \rangle \, dt \nonumber \\[1ex]
&=& 2 \int^{\infty}_{- \infty} \int^{l_y}_0 \int^{l_x}_0 
    R_e \left\{(\psi_2 \hat{\psi}_2)  f + (\psi_1 \hat{\psi}_1) 
      \overline{f} \right\} \, dx \, dy \, dt \, ,  \label{mlf2} 
\end{eqnarray}
where the integrand is evaluated on an unperturbed heteroclinic
orbit obtained through the B\"acklund-Darboux
transformations given in Theorem \ref{DSTH}. A concrete example 
has been worked out in \cite{Li00a}.

\section{Two-Dimensional Euler Equations}

One of the most important implications of chaos theory of partial 
differential equations in theoretical physics will be on the study 
of turbulence. For that goal, the author choose the 2D Navier-Stokes 
equations under periodic boundary conditions to begin a dynamical system 
study.
\begin{eqnarray}
 & &{\pa \Om \over \pa t}= - u \ {\pa \Om \over \pa x} -v
\ {\pa \Om \over \pa y} +\e \bigg [ \Delta \Om + f \bigg ],
\label{NSE} \\
& &{\pa u \over \pa x} +{\pa v \over \pa y} = 0; \nonumber
\end{eqnarray}
under periodic boundary conditions in both $x$ and $y$ directions
with period $2\pi$, where $\Om$ is vorticity, $u$ and $v$ are 
respectively velocity components along $x$ and $y$ directions, $\e = 1/Re$,
and $f$ is the body force.
When $\e =0$, we have the 2D Euler equations, 
\begin{eqnarray}
& &{\pa \Om \over \pa t}= - u \ {\pa \Om \over \pa x} -v
\ {\pa \Om \over \pa y}\ ,
\label{Euler} \\
& &{\pa u \over \pa x} +{\pa v \over \pa y} = 0\ . \nonumber
\end{eqnarray}
The relation between vorticity $\Om$ and stream 
function $\Psi$ is,
\[
\Om ={\pa v \over \pa x} -{\pa u \over \pa y} =\Dl \Psi \ ,
\]
where the stream function $\Psi$ is defined by,
\[
u=- {\pa \Psi \over \pa y} \ ,\ \ \ v={\pa \Psi \over \pa x} \ .
\]

\subsection{Lax Pair and Darboux Transformation}

The main breakthrough in this project is the discovery of the Lax pair 
for 2D Euler equation \cite{Li01a}. The philosophical 
significance of the existence of a 
Lax pair for 2D Euler equation is beyond the particular 
project undertaken here. If one defines integrability of an equation 
by the existence of a Lax pair, then 2D Euler equation 
is integrable. More importantly, 2D Navier-Stokes equation 
at high Reynolds numbers is a near integrable system. Such a point of view 
changes our old ideology on Euler and Navier-Stokes equations.

Starting from Lax pairs, homoclinic structures can be constructed through 
Darboux transformations \cite{Li00a}. Indeed, in \cite{LY00}, the Darboux 
transformation 
for the Lax pair of 2D Euler equation has been found. Our general program 
is to first identify the figure eight structures of 2D Euler equation, and 
then study their consequence in 2D Navier-Stokes equation. 
The high Reynolds number 2D Navier-Stokes 
equation is viewed as a singular perturbation of the 2D Euler equation through 
the perturbation $\e \Dl$, where $\e = 1/Re$ is the inverse of the Reynolds 
number. As mentioned above, singular perturbations have been investigated 
for nonlinear Schr\"odinger equations.

We consider the 2D Euler equation,
\begin{equation}
{\pa \Om \over \pa t} + \{ \Psi, \Om \} = 0 \ ,
\label{euler}
\end{equation}
where the bracket $\{\ ,\  \}$ is defined as
\[
\{ f, g\} = (\pa_x f) (\pa_y g) - (\pa_y f) (\pa_x g) \ , \ \ \mbox{and}\
\Om = \Dl \Psi \ .
\]
\begin{theorem}[\cite{Li01a}]
The Lax pair of the 2D Euler equation (\ref{euler}) is given as
\begin{equation}
\left \{ \begin{array}{l} 
L \varphi = \la \varphi \ ,
\\
\pa_t \varphi + A \varphi = 0 \ ,
\end{array} \right.
\label{laxpair}
\end{equation}
where
\[
L \varphi = \{ \Om, \varphi \}\ , \ \ \ A \varphi = \{ \Psi, \varphi \}\ ,
\]
and $\la$ is a complex constant, and $\varphi$ is a complex-valued function.
\end{theorem}
In \cite{LY00}, A 
B\"acklund-Darboux transformation is found for the above Lax pair.
Consider the Lax pair (\ref{laxpair}) at $\la =0$, i.e.
\begin{eqnarray}
& & \{ \Om, p \} = 0 \ , \label{d1} \\
& & \pa_t p + \{ \Psi, p \} = 0 \ , \label{d2} 
\end{eqnarray}
where we replaced the notation $\varphi$ by $p$.
\begin{theorem}
Let $f = f(t,x,y)$ be any fixed solution to the system 
(\ref{d1}, \ref{d2}), we define the Gauge transform $G_f$:
\begin{equation}
\tilde{p} = G_f p = \frac {1}{\Om_x}[p_x -(\pa_x \ln f)p]\ ,
\label{gauge}
\end{equation}
and the transforms of the potentials $\Om$ and $\Psi$:
\begin{equation}
\tilde{\Psi} = \Psi + F\ , \ \ \ \tilde{\Om} = \Om + \Dl F \ ,
\label{ptl}
\end{equation}
where $F$ is subject to the constraints
\begin{equation}
\{ \Om, \Dl F \} = 0 \ , \ \ \ \{ \Om, F \} = 0\ .
\label{constraint}
\end{equation}
Then $\tilde{p}$ solves the system (\ref{d1}, \ref{d2}) at 
$(\tilde{\Om}, \tilde{\Psi})$. Thus (\ref{gauge}) and 
(\ref{ptl}) form the Darboux transformation for the 2D 
Euler equation (\ref{euler}) and its Lax pair (\ref{d1}, \ref{d2}).
\label{dt}
\end{theorem}

\subsection{Linearized 2D Euler Equations}

Under the periodic boundary condition and requiring that both 
$u$ and $v$ have means zero,
\[
\int_0^{2\pi}\int_0^{2\pi} u\ dxdy =\int_0^{2\pi}\int_0^{2\pi} v\ dxdy=0,
\]
expanding $\Om$ into Fourier series,
$\Om =\sum_{k\in Z^2/\{0\}} \om_k \ e^{ik\cdot X}$,
where $\om_{-k}=\overline{\om_k}\ $, $k=(k_1,k_2)$, 
$X=(x,y)$, the system (\ref{Euler}) can be rewritten as the 
following kinetic system,
\begin{equation}
\dot{\om}_k = \sum_{k=p+q} A(p,q) \ \om_p \om_q \ ,
\label{Keuler}
\end{equation}
where $A(p,q)$ is given by,
\begin{eqnarray}
A(p,q)&=& {1\over 2}[|q|^{-2}-|p|^{-2}](p_1 q_2 -p_2 q_1) \label{Af} \\      
      &=& {1\over 2}[|q|^{-2}-|p|^{-2}]\left | \begin{array}{lr} 
p_1 & q_1 \\ p_2 & q_2 \\ \end{array} \right | \ , \nonumber
\end{eqnarray}
where $|q|^2 =q_1^2 +q_2^2$ for $q=(q_1,q_2)$, similarly for $p$.
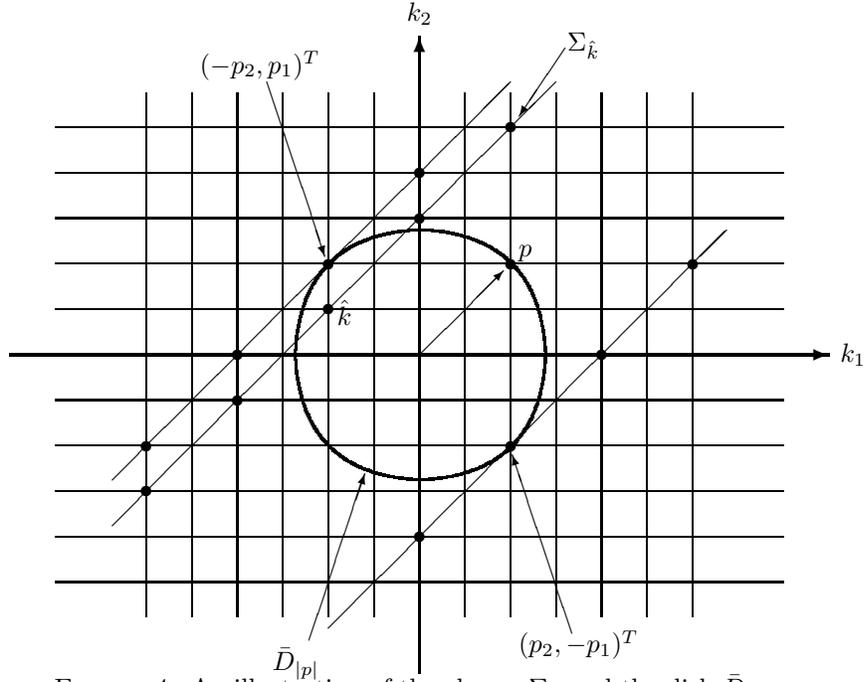
\begin{figure}[ht]
  \begin{center}
    \leavevmode
      \setlength{\unitlength}{2ex}
  \begin{picture}(36,27.8)(-18,-12)
    \thinlines
\multiput(-12,-11.5)(2,0){13}{\line(0,1){23}}
\multiput(-16,-10)(0,2){11}{\line(1,0){32}}
    \thicklines
\put(0,-14){\vector(0,1){28}}
\put(-18,0){\vector(1,0){36}}
\put(0,15){\makebox(0,0){$k_2$}}
\put(18.5,0){\makebox(0,0)[l]{$k_1$}}
\qbezier(-5.5,0)(-5.275,5.275)(0,5.5)
\qbezier(0,5.5)(5.275,5.275)(5.5,0)
\qbezier(5.5,0)(5.275,-5.275)(0,-5.5)
\qbezier(0,-5.5)(-5.275,-5.275)(-5.5,0)
    \thinlines
\put(4,4){\circle*{0.5}}
\put(0,0){\vector(1,1){3.7}}
\put(4.35,4.35){$p$}
\put(4,-4){\circle*{0.5}}
\put(8,0){\circle*{0.5}}
\put(-8,0){\circle*{0.5}}
\put(-8,-2){\circle*{0.5}}
\put(-12,-4){\circle*{0.5}}
\put(-12,-6){\circle*{0.5}}
\put(-4,2){\circle*{0.5}}
\put(-4,4){\circle*{0.5}}
\put(0,6){\circle*{0.5}}
\put(0,8){\circle*{0.5}}
\put(4,10){\circle*{0.5}}
\put(12,4){\circle*{0.5}}
\put(0,-8){\circle*{0.5}}
\put(-4,-12){\line(1,1){17.5}}
\put(-13.5,-7.5){\line(1,1){19.5}}
\put(-13.5,-5.5){\line(1,1){17.5}}
\put(-3.6,1.3){$\hat{k}$}
\put(-7,12.1){\makebox(0,0)[b]{$(-p_2, p_1)^T$}}
%\put(-7.85,11.75){\vector(1,-2){3.65}}
\put(-6.7,12){\vector(1,-3){2.55}}
\put(6.5,13.6){\makebox(0,0)[l]{$\Sg_{\hat{k}}$}}
\put(6.4,13.5){\vector(-2,-3){2.0}}
\put(7,-12.1){\makebox(0,0)[t]{$(p_2, -p_1)^T$}}
\put(6.7,-12.25){\vector(-1,3){2.62}}
\put(-4.4,-13.6){\makebox(0,0)[r]{$\bar{D}_{|p|}$}}
\put(-4.85,-12.55){\vector(1,3){2.45}}
  \end{picture}
  \end{center}
\caption{An illustration of the classes $\Sg_{\hk}$ and the disk 
$\bar{D}_{|p|}$.}
\label{class}
\end{figure}
\begin{figure}[ht]
  \begin{center}
    \leavevmode
      \setlength{\unitlength}{2ex}
  \begin{picture}(36,27.8)(-18,-12)
    \thicklines
\put(0,-14){\vector(0,1){28}}
\put(-18,0){\vector(1,0){36}}
\put(0,15){\makebox(0,0){$\Im \{ \la \}$}}
\put(18.5,0){\makebox(0,0)[l]{$\Re \{ \la \}$}}
\put(0.1,-7){\line(0,1){14}}
\put(.2,-.2){\makebox(0,0)[tl]{$0$}}
\put(-0.2,-7){\line(1,0){0.4}}
\put(-0.2,7){\line(1,0){0.4}}
\put(2.0,-6.4){\makebox(0,0)[t]{$-i2|b|$}}
\put(2.0,7.6){\makebox(0,0)[t]{$i2|b|$}}
\end{picture}
  \end{center}
\caption{The spectrum of $\LL_A$ in case 1.}
\label{splb}
\end{figure}
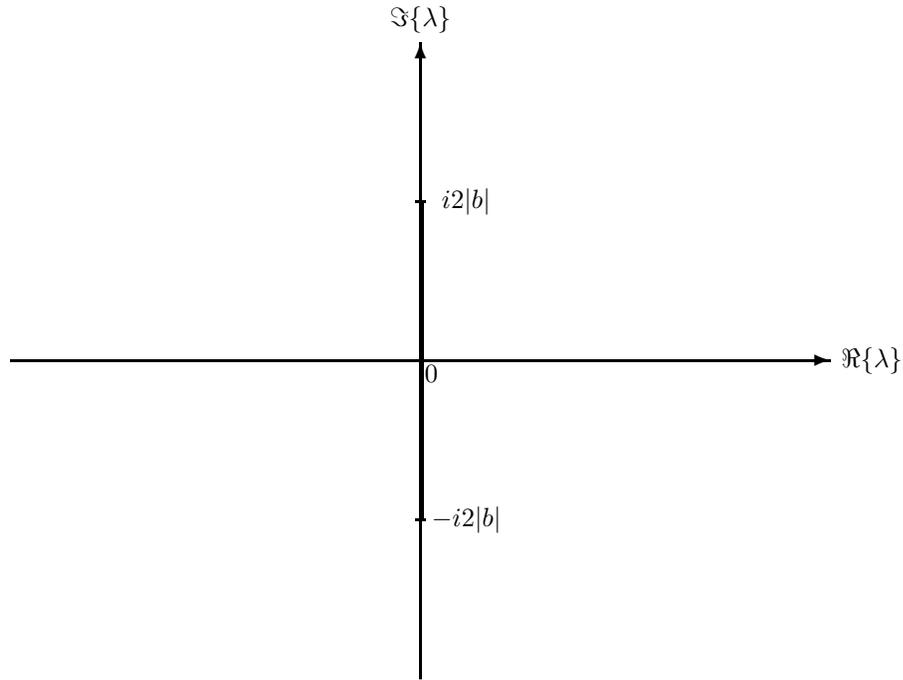
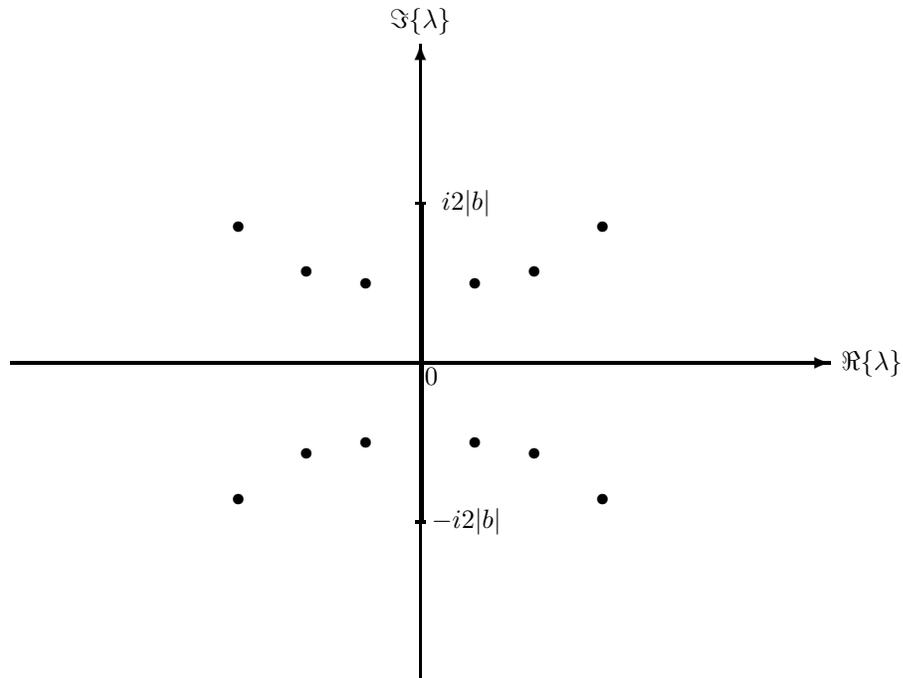
\begin{figure}[ht]
  \begin{center}
    \leavevmode
      \setlength{\unitlength}{2ex}
  \begin{picture}(36,27.8)(-18,-12)
    \thicklines
\put(0,-14){\vector(0,1){28}}
\put(-18,0){\vector(1,0){36}}
\put(0,15){\makebox(0,0){$\Im \{ \la \}$}}
\put(18.5,0){\makebox(0,0)[l]{$\Re \{ \la \}$}}
\put(0.1,-7){\line(0,1){14}}
\put(.2,-.2){\makebox(0,0)[tl]{$0$}}
\put(-0.2,-7){\line(1,0){0.4}}
\put(-0.2,7){\line(1,0){0.4}}
\put(2.0,-6.4){\makebox(0,0)[t]{$-i2|b|$}}
\put(2.0,7.6){\makebox(0,0)[t]{$i2|b|$}}
\put(2.4,3.5){\circle*{0.5}}
\put(-2.4,3.5){\circle*{0.5}}
\put(2.4,-3.5){\circle*{0.5}}
\put(-2.4,-3.5){\circle*{0.5}}
\put(5,4){\circle*{0.5}}
\put(-5,4){\circle*{0.5}}
\put(5,-4){\circle*{0.5}}
\put(-5,-4){\circle*{0.5}}
\put(8,6){\circle*{0.5}}
\put(-8,6){\circle*{0.5}}
\put(8,-6){\circle*{0.5}}
\put(-8,-6){\circle*{0.5}}
\end{picture}
  \end{center}
\caption{The spectrum of $\LL_A$ in case 2.}
\label{spla2}
\end{figure}
\begin{figure}[ht]
  \begin{center}
    \leavevmode
      \setlength{\unitlength}{2ex}
  \begin{picture}(36,27.8)(-18,-12)
    \thicklines
\put(0,-14){\vector(0,1){28}}
\put(-18,0){\vector(1,0){36}}
\put(0,15){\makebox(0,0){$\Im \{ \tla \}$}}
\put(18.5,0){\makebox(0,0)[l]{$\Re \{ \tla \}$}}
\put(2.4,3.5){\circle*{0.5}}
\put(-2.4,3.5){\circle*{0.5}}
\put(2.4,-3.5){\circle*{0.5}}
\put(-2.4,-3.5){\circle*{0.5}}  
\end{picture}
  \end{center}
\caption{The quadruple of eigenvalues 
for the system led by the class $\Sg_{\hat{k}}$ labeled by $\hk = (1,0)^T$, 
when $p=(1,1)^T$.}
\label{figev}
\end{figure}
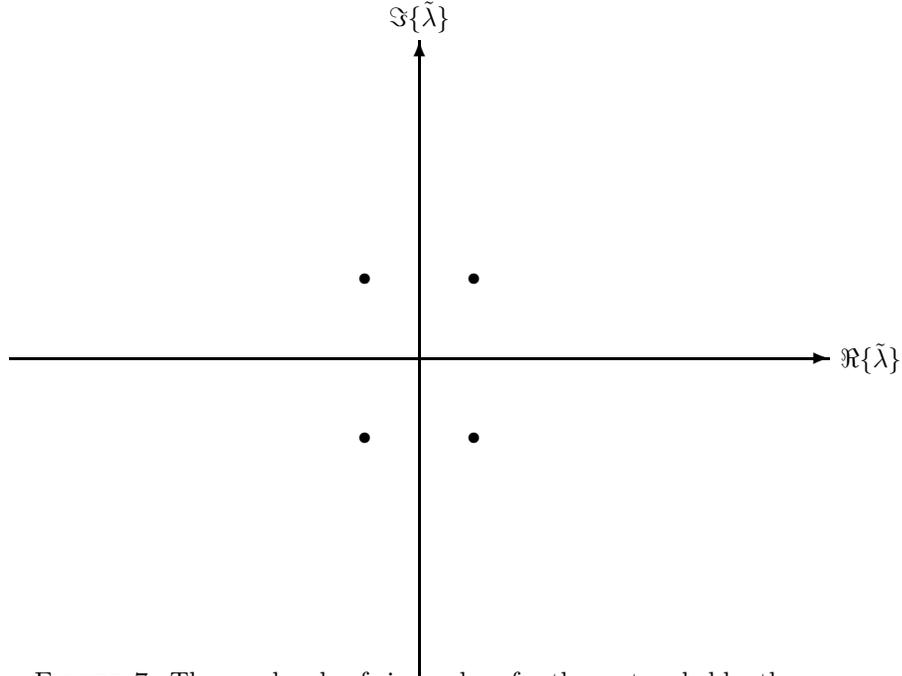
\begin{figure}
\vspace{7.5in}
\caption{The heteroclinic orbits and unstable manifolds of the 
Galerkin truncation.}
\label{apr}
\end{figure}
To understand the hyperbolic structures of the 2D Euler equations, we 
first investigate the linearized 2D Euler equations at a stationary solution. 
Denote $\{ \om_k \}_{k\in \Z}$ by $\om$. Consider the simple fixed point 
$\om^*$:
\begin{equation}
\om^*_p = \Ga,\ \ \ \om^*_k = 0 ,\ \mbox{if} \ k \neq p \ \mbox{or}\ -p,
\label{fixpt}
\end{equation}
of the 2D Euler equation (\ref{Keuler}), where 
$\Ga$ is an arbitrary complex constant. The 
{\em{linearized two-dimensional Euler equation}} at $\om^*$ is given by,
\begin{equation}
\dot{\om}_k = A(p,k-p)\ \Ga \ \om_{k-p} + A(-p,k+p)\ \bar{\Ga}\ \om_{k+p}\ .
\label{LE}
\end{equation}
\begin{definition}[Classes]
For any $\hk \in \Z$, we define the class $\Sg_{\hk}$ to be the subset of 
$\Z$:
\[
\Sg_{\hk} = \bigg \{ \hk + n p \in \Z \ \bigg | \ n \in Z, \ \ p \ \mbox{is 
specified in (\ref{fixpt})} \bigg \}.
\]
\label{classify}
\end{definition}
See Figure \ref{class} for an illustration of the classes. 
According to the classification 
defined in Definition \ref{classify}, the linearized two-dimensional Euler 
equation (\ref{LE}) decouples into infinite many invariant subsystems:
\begin{eqnarray}
\dot{\omega}_{\hat{k} + np} &=& A(p, \hat{k} + (n-1) p) 
     \ \Gamma \ \omega_{\hat{k} + (n-1) p} \nonumber \\  
& & + \ A(-p, \hat{k} + (n+1)p)\ 
     \bar{\Gamma} \ \omega_{\hat{k} +(n+1)p}\ . \label{CLE}
\end{eqnarray}
\begin{definition}[The Disk]
The disk of radius $\left| p \right|$ in $Z^2 / \left\{ 0
\right\}$, denoted by $D_{\left| p \right|}$, is defined as
\[
  D_{\left| p \right|} = \bigg \{ k \in Z^2 / \left\{ 0 \right\} \ \bigg| 
      \ \left| k \right| < \left| p \right| \bigg \} \, .
\]
The closure of $D_{\left| p \right|}$, denoted by
$\bar{D}_{\left| p \right|}$, is defined as
\[ 
 \bar{D}_{\left| p \right|} = \bigg \{ k \in Z^2/ \left\{ 0 \right\} \ \bigg| 
     \ \left| k \right| \leq \left| p \right| \bigg \} \, .
\]
\end{definition}
\begin{theorem}[Unstable Disk Theorem]
If $\ \Sigma_{\hat{k}} \cap \bar{D}_{\left| p \right|} = \emptyset,\ $
then the invariant subsystem (\ref{CLE}) is Liapunov stable
for all $t \in R$, in fact, 
\[
  \sum_{n \in Z} \left|
  \omega_{\hat{k}+np}(t) \right|^2 \leq \sigma \ \sum_{n \in Z} \left|
  \omega_{\hat{k}+np}(0) \right|^2 \, , \quad \quad \forall t \in R \, , 
\]
where 
\[
  \sigma = \left[ \max_{n \in Z} 
              \left\{ - \rho_n \right\}
           \right] \, 
           \left[ \min_{n \in Z} 
               \left\{ -\rho_n \right\}
           \right]^{-1} \, , \quad 0< \sigma < \infty \, .
\]
\label{UDT}
\end{theorem}
\begin{theorem}
The eigenvalues of the linear system (\ref{CLE}) are of 
four types: real pairs ($c, -c$), purely imaginary pairs ($id, -id$), 
quadruples ($\pm c \pm id$), and zero eigenvalues.
\end{theorem}
\begin{theorem}[The Spectral Theorem]
\begin{enumerate}
\item If $\Sg_{\hat{k}} \cap \bar{D}_{|p|} = \emptyset$, then the entire
$\ell_2$ spectrum of the linear operator $\LL_A$ (defined by the right-hand 
side of the invariant subsystem) is its continuous spectrum. See Figure 
\ref{splb}.
\item If $\Sg_{\hat{k}} \cap \bar{D}_{|p|} \neq \emptyset$, then the entire
essential $\ell_2$ spectrum of the linear operator $\LL_A$ is its 
continuous spectrum. That is, the residual spectrum of $\LL_A$ is empty, 
$\sg_r (\LL_A) = \emptyset$. The point spectrum 
of $\LL_A$ is symmetric with respect to both real and imaginary axes. 
See Figure \ref{spla2}.
\end{enumerate}
\label{spthla}
\end{theorem}
We can calculate the eigenvalues through continued fractions.
Let $p=(1,1)^T$, in this case, only one class $\Sg_{\hat{k}}$ 
labeled by $\hk = (1,0)^T$ has no empty intersection with $\bar{D}_{|p|}$ 
(the other class labeled by $\hk = (0,1)^T$ gives the complex conjugate 
of the system led by the class labeled by $\hk = (1,0)^T$). For this class,
there is no real eigenvalue. Numerical 
calculation through continued fractions gives the eigenvalue:
\[
\tla=0.24822302478255 \ + \ i \ 0.35172076526520\ .
\]
Thus we have a quadruple 
of eigenvalues, see Figure \ref{figev} for an illustration.
Denote by $L$ the right hand side of (\ref{LE}), the spectral mapping theorem 
holds.
\begin{theorem}[\cite{LLM01}]
$$\sigma(e^{tL})=e^{t\sigma(L)}, t\neq 0.$$
\end{theorem}
Moreover, the number of eigenvalues has a sharp upper bound.
Let $\z$ denote the number of points $q \in
\Z$ that belong to the open disk of radius $|p|$, 
and such that $q$ is not parallel to $p$. 
\begin{theorem}[\cite{LLM01}]
The number of nonimaginary eigenvalues of $L$ (counting the multiplicities) 
does not exceed $2\z$.
\end{theorem}

\subsection{Approximate Explicit Representations of the Hyperbolic Structures
of 2D Euler Equations}

From Figure \ref{figev}, we see that the simple fixed point given by 
$p = (1,1)$, has unstable eigenvalues. Our interest is to obtain 
representations
of the correponding hyperbolic structures for 2D Euler equations. 
In \cite{Li01d}, through Galerkin truncation, we obtained the approximate 
explicit representation. Figure \ref{apr} shows the heteroclinic 
orbits and unstable manifolds of the Galerkin truncation.

\section{Conclusion and Discussion}

We have reported the status of chaos in nonlinear wave 
equations and of study on 2D Euler equations. In particular, 
we have summarized the most recent results on Lax pair and Darboux 
transformations for 2D Euler equations.


\begin{thebibliography}{99}

\bibitem{AC91} M.J. Ablowitz and P. A. Clarkson, {\it Solitons, Nonlinear 
Evolution Equations and Inverse Scattering}, London Mathematical Society 
Lecture Note Series 149, 1991.
\bibitem{AL76} M.J. Ablowitz and J. F. Ladik, A nonlinear difference scheme 
and inverse scattering, {\it Stud. Appl. Math.} {\bf 55} (1976), 213.
\bibitem{AS81} M.J. Ablowitz and H. Segur, {\it Solitons and the 
Inverse Scattering Transform}, SIAM, Philadelphia, 1981.
\bibitem{Bou94} J. Bourgain, Construction of quasi-periodic solutions for
		  Hamiltonian perturbations of linear
		  equations and applications to nonlinear PDE,
{\it International Mathematics Research Notes} {\bf 11} (1994), 1.
\bibitem{CW93} W. Craig and C. E. Wayne, Newton's method and periodic 
solutions of nonlinear wave equation, {\it Comm. Pure Appl. Math.} 
{\bf 46} (1993), 1409.
\bibitem{DS74} A. Davey and K. Stewartson, On three-dimensional packets of 
surface waves, {\it Proc. R. Soc. Lond.} {\bf A338} (1974), 101.
\bibitem{Gra74} S. M. Graff, On the construction of hyperbolic invariant
		  tori for Hamiltonian systems, {\it J. Diff. Eqs.} 
{\bf 15} (1974), 1.
\bibitem{GH83} J. Guckenheimer and P. J. Holmes, {\it Nonlinear 
Oscillations, Dynamical Systems, and Bifurcations of Vector Fields},
Springer-Verlag, Applied Mathematical Sciences, vol.42, 1983.
\bibitem{Kuk93} S. B. Kuksin, {\it Nearly Integrable Infinite-Dimensional
		  Hamiltonian Systems}, Springer-Verlag, Lecture Notes 
Math. 1556, 1993.
\bibitem{LLM01} Y. Latushkin and Y. Li and M. Stanislavova, The spectrum 
of a linearized 2D Euler operator, {\it Submitted} (2001).
\bibitem{Li92} Y. Li, Backlund transformations and homoclinic 
structures for the NLS equation, {\it Phys. Letters A} {\bf 163}
(1992), 181.
\bibitem{Li98a} Y. Li, Homoclinic tubes and symbolic dynamics 
           discrete nonlinear Schr{\"{o}}dinger equation 
           under Hamiltonian perturbations, {\it Submitted} (1998).
\bibitem{Li99b} Y. Li, Homoclinic tubes in nonlinear Schr{\"{o}}dinger 
		  equation under Hamiltonian perturbations, {\it 
Progress of Theoretical Physics} {\bf 101, No. 3(4)} (1999), 559.
\bibitem{Li99a} Y. Li, Smale horseshoes and symbolic dynamics 
          in perturbed nonlinear Schr{\"{o}}dinger equations, {\it 
Journal of Nonlinear Sciences} {\bf 9, no.4} (1999), 363.
\bibitem{Li00a} Y. Li, B{\"{a}}cklund-Darboux transformations 
         and Melnikov analysis for Davey-Stewartson II equations, {\it
Journal of Nonlinear Sciences} {\bf 10, no.1} (2000), 103.
\bibitem{Li00b} Y. Li, On 2D Euler equations. I. On the energy-Casimir 
           stabilities and the spectra for a linearized two 
           dimensional Euler equation, {\it Journal of Mathematical 
Physics} {\bf 41, no.2} (2000), 728.
\bibitem{Li01a} Y. Li, A Lax pair for the 2D Euler equation, {\it 
Journal of Mathematical Physics} {\bf 42, no.8} (2001), 3552.
\bibitem{Li01d} Y. Li, On 2D Euler equations: part II. Lax pairs 
and homoclinic strucutres, {\it Submitted} (2001). 
\bibitem{Li01b} Y. Li, Persistent homoclinic orbits for nonlinear 
Schr{\"{o}}dinger equation under singular perturbation, {\it Submitted} 
(2001).
\bibitem{Li01c} Y. Li, Singularly perturbed vector and scalar nonlinear 
Schr{\"{o}}dinger equation with persistent homoclinic orbits, {\it 
Studies in Applied Mathematics, to appear} (2001).
\bibitem{LMSW96} Y. Li et al., Persistent homoclinic orbits for a 
perturbed nonlinear Schr{\"{o}}dinger equation, 
{\it Comm. Pure Appl. Math.} {\bf XLIX} (1996), 1175.
\bibitem{LM94} Y. Li and D. W. McLaughlin, Morse and Melnikov functions 
for NLS pdes, {\it Comm. Math. Phys.} {\bf 162} (1994), 175.
\bibitem{LM97} Y. Li and D. W. McLaughlin, Homoclinic orbits and chaos 
in discretized perturbed NLS system, part I. homoclinic orbits,
{\it Journal of Nonlinear Sciences} {\bf 7} (1997), 211.
\bibitem{LW97a} Y. Li and S. Wiggins, Homoclinic orbits and chaos 
in discretized perturbed NLS system, part II. symbolic dynamics,
{\it Journal of Nonlinear Sciences} {\bf 7} (1997), 315.
\bibitem{LY00} Y. Li and A. Yurov, Lax pairs and Darboux transformations 
for Euler equations, {\bf Preprint, available at: 
http://xxx.lanl.gov/abs/math.AP/0101214, or 
http://www.math.missouri.edu/\~{}cli} (2000).
\bibitem{LL83} A. J. Lichtenberg and M. A. Lieberman, {\it 
Regular and Stochastic Motion}, Springer-Verlag, 
Applied Mathematical Sciences, vol.38, 1983.
\bibitem{Oza92} T. Ozawa, Exact blow-up solutions to the Cauchy
		  problem for the Davey-Stewartson systems,
{\it Proc. R. Soc. Lond.} {\bf 436} (1992), 345.
\bibitem{Pos96} J. Poschel, A KAM-theorem for some nonlinear PDEs, 
{\it Ann. Scuola Norm. Sup. Pisa, Cl. Sci., IV Ser. 15} {\bf
23} (1996), 119.
\bibitem{Way90} C. E. Wayne, Periodic and quasi-periodic solutions of
		  nonlinear wave equations via KAM theory, 
{\it Commun. Math. Phys.} {\bf 127} (1990), 479.
\bibitem{Zen00b} C. Zeng, Erratum: Homoclinic orbits for a perturbed 
nonlinear Schr{\"{o}}dinger equation, {\it Comm. Pure Appl. Math.} {\bf 
53, no.12} (2000), 1603.
\bibitem{Zen00a} C. Zeng, Homoclinic orbits for a perturbed 
nonlinear Schr{\"{o}}dinger equation, {\it Comm. Pure Appl. Math.} {\bf 
53, no.10} (2000), 1222.











\end{thebibliography}
\end{document}